\documentclass[11pt]{amsart}
\usepackage{amssymb}
\usepackage{amsmath}
\usepackage{amscd}
\usepackage{curves}
\usepackage{epsfig}
\usepackage{geometry}
\usepackage{graphicx}
\usepackage{pstricks}
\usepackage{pstricks-add}
\usepackage{pst-grad}
\usepackage{pst-plot}
\usepackage{verbatim}
\usepackage{latexsym,eucal,amsfonts,amssymb,amsmath,graphicx}
\numberwithin{equation}{section}
\newtheorem{theorem}{Theorem}[section]
\newtheorem{lemma}[theorem]{Lemma}
\newtheorem{prop}[theorem]{Proposition}

\def \bpf {\begin{proof}}
\def \epf {\end{proof}}
\def \beq {\begin{equation*}}
\def \eeq {\end{equation*}}
\def \bsp{\begin{split}}
\def \esp{\end{split}}

\def \mcf {{\mathcal F}}
\def \mcg {{\mathcal G}}

\def \mcm {{\mathcal M}}

\def \txi {\widetilde{\xi}}

\def \mh {{\mathbb H}}
\def \mb {{\mathbb B}}

\def \mbr {{\mathbb R}}
\def \mr {{\mathbb R}}

\def\ha {\frac{1}{2}}

\def \intx {\mathring{X}}

\def \ff {\text{ff}}

\def \diag{\textrm{Diag}}
\def \odiag{\overline{\textrm{Diag}}}

\def \mrn {{\mathbb R}^n}

\def \eps {\varepsilon}   
\def \la {\lambda}   
\def \La {\Lambda}

\def \del {\delta}   

\def \p {\partial}
\def \phan {\phantom{a}}

\def \xo {X\times_0 X}

\numberwithin{equation}{section}

\begin{document}
\title[The Scattering Relation on  AHM]{The Scattering Relation on Asymptotically Hyperbolic Manifolds}
\author{Ant\^onio S\'a Barreto and Yiran Wang}
\email{sabarre@math.purdue.edu, wang554@purdue.edu}
\address{Department of Mathematics, Purdue University \newline
\indent 150 North University Street, West Lafayette Indiana,  47907, USA}
\begin{abstract}  We study the scattering relation and the sojourn times on non-trapping asymptotically hyperbolic manifolds and use it  to obtain the asymptotics of the distance function on  geodesically convex asymptotically hyperbolic manifolds.    
\end{abstract}
\maketitle

\section{Introduction}

 Asymptotically hyperbolic manifolds are a generalization of the  hyperbolic space $(\mb^{n+1},g),$  
\begin{equation}\label{hyper}
\bsp
&\mb^{n+1}=\{ z\in \mbr^{n+1}: |z|<1\} \text{ equipped with the metric } g=\frac{4dz^2}{(1-|z|^2)^2}.
\end{split}
\end{equation}
 Instead of $\mb^{n+1},$ we consider the interior of a  $C^\infty$ manifold with boundary $X$ of dimension $n+1,$ and assume that the interior of $X,$ which we denote by $\intx,$ is equipped with a metric $g$ such that for any defining function $\rho$ of  $\p X$  (i.e. $\rho\in C^\infty (X),$  $\rho>0$ in  $\intx,$  $\{\rho=0\}= \p X$ and $d\rho\neq 0$ at $\partial X$),  $\rho^2g$ is a $C^\infty$  non-degenerate Riemannian metric up to $\p X.$ In the case of the hyperbolic space $X={\overline{\mb}}^{n+1}$ and $\rho=1-|z|^2.$    According to \cite{MM} the manifold $(\intx,g)$ is  complete and its sectional curvatures approach $-\left| d\rho|_{\p X}\right|_{h_0}^2,$ as $\rho\downarrow 0$ along any curve, where  $h_0=\rho^2 g|_{\p X}.$    In particular, when 
\begin{gather}
\displaystyle \left|d\rho|_{\p X}\right|_{h_0}=1, \label{asympt-curv}
\end{gather}  
the sectional curvature converges to $-1$ at the boundary.    Following Mazzeo and Melrose \cite{MM},  manifolds $(\intx,g)$ for which $\rho^2g$ is non-degenerate at $\p X$ and  \eqref{asympt-curv} holds 
are called asymptotically hyperbolic manifolds (AHM).   It follows from the definition that the metric $g$  determines a conformal structure on $\p X,$ and because of that these manifolds have been studied in connection with conformal field theory \cite{Gr,GrZw}.  As shown in \cite{Gr,JS1},  if $h_0\in [\rho^2g|_{\p X}],$ the equivalence class of $\rho^2g|_{\p X},$  there exists a boundary defining function $x$ in a neighborhood of $\p X$ such that
\begin{equation}\label{prod}
g = \frac{dx^2}{x^2} +\frac{ h(x)}{x^2},  \;\ h(0)=h_0, \text{ on }   [0,\eps) \times \p X,
\end{equation}
where $h(x)$ is a $C^\infty$ family of Riemannian metrics on $\p X$ parametrized by $x.$    Of course, $x$ can be extended (non-uniquely) to $X$ by setting it equal to a constant on a compact set of $\intx.$

We say that $(\intx,g)$ is non-trapping if any geodesic $\gamma(t) \rightarrow \p X$ as $\pm t\rightarrow \infty,$ and  we shall assume  throughout this paper that $(\intx, g)$ is non-trapping.  
Our goal is to understand the behavior of geodesics on non-trapping AHM and define the scattering relation at the boundary at infinity.  
One can easily describe the scattering relation for non-trapping compactly supported metric perturbations of Euclidean space.  Suppose that $g=\sum_{i,j=1}^n g_{ij}(x) dx_idx_i$ is a $C^\infty$ non-trapping Riemannian metric on $\mrn$ and suppose that $g_{ij}(x)=\del_{ij}$ if $x\not\in K\subset \mrn,$ where $K$ is compact.  Let $B$ be a bounded ball and suppose $K\subset B.$   A light ray comes from $\mrn\setminus B$ enters $B,$ is scattered by the metric in $K$ and comes out of $B.$ If the light ray  that comes into $B$  intersects $\p B$ at a point $z\in \p B$ in the direction $\zeta$ and comes out at points $z'\in \p B$ with direction $\zeta',$ the map  $(z,\zeta)\longmapsto (z',\zeta')$  is called the scattering relation.  The time $t$ that it takes for the geodesic to travel across $B$ is called the sojourn time.  This can also be described in terms of  the submanifold 
$\Lambda=\{ (t,1,z,\zeta, z',\zeta'): (z,\zeta)=\exp(t H_q)(z',\zeta')\},$ where $q=\sum_{i,j=1}^n g^{ij}\xi_i\xi_j,$ $g^{-1}=(g^{ij})$ is the dual metric to $g.$  The scattering relation is then the restriction of $\Lambda$ to $\p B \times \p B,$ i.e $\Lambda\cap(T_{\p B}^* \mrn \times T_{\p B}^* \mrn).$  Theorem \ref{soj}  below can be used to define the analogue of such a map for  non-trapping AHM.    Let $\Lambda_1$  be the extension of 
$\beta_1^*\tilde{\Lambda}$ up to $\p(\mr_s \times \xo),$ as defined below.  The scattering relation is 
$\Lambda_1 \cap T_{\mr_s\times \{\rho_R=\rho_L=0\}}^* (\mr_s\times \xo).$  
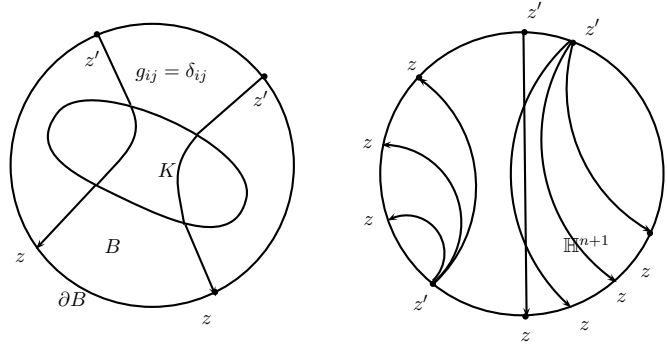
\begin{figure}
\scalebox{.7} 
{
\begin{pspicture}(0,-3.6123857)(12.461016,2.9534347)
\psbezier[linewidth=0.04](0.9410156,0.8150362)(0.32046875,0.03086683)(1.3198191,-0.41155323)(2.2210157,-0.8449638)(3.1222122,-1.2783744)(4.3186464,-1.7842731)(4.5210156,-0.8049638)(4.723385,0.1743455)(1.5615625,1.5992055)(0.9410156,0.8150362)
\pscircle[linewidth=0.04,dimen=outer](2.7310157,-0.19496381){2.71}
\psdots[dotsize=0.12](4.861016,1.4950362)
\usefont{T1}{ptm}{m}{n}
\rput(4.792471,1.0800362){$z'$}
\psdots[dotsize=0.12](1.6810156,2.295036)
\usefont{T1}{ptm}{m}{n}
\rput(1.5924706,1.8400362){$z'$}
\psdots[dotsize=0.12](3.9210157,-2.6049638)
\usefont{T1}{ptm}{m}{n}
\rput(3.7624707,-3.119964){$z$}
\usefont{T1}{ptm}{m}{n}
\rput(3.0124707,-0.27996382){$K$}
\usefont{T1}{ptm}{m}{n}
\rput(1.2124707,-2.7399638){$\p B$}
\usefont{T1}{ptm}{m}{n}
\rput(1.9824708,-1.7199638){$B$}
\pscircle[linewidth=0.04,dimen=outer](9.751016,-0.3549638){2.71}
\psline[linewidth=0.04cm,arrowsize=0.05291667cm 2.0,arrowlength=1.4,arrowinset=0.4]{->}(9.781015,2.3350363)(9.841016,-3.0849638)
\rput{-46.93517}(2.465871,4.469973){\psarc[linewidth=0.04,arrowsize=0.05291667cm 2.0,arrowlength=1.4,arrowinset=0.4]{->}(6.381016,-0.60496384){2.5}{0.0}{102.991585}}
\usefont{T1}{ptm}{m}{n}
\rput(7.832471,-2.7999637){$z'$}
\psdots[dotsize=0.12](7.801016,1.4750361)
\psdots[dotsize=0.12](8.061016,-2.444964)
\psdots[dotsize=0.12](9.801016,2.3350363)
\psdots[dotsize=0.12](10.721016,2.1350362)
\psdots[dotsize=0.12](12.181016,-1.4849638)
\usefont{T1}{ptm}{m}{n}
\rput(7.662471,1.7000362){$z$}
\psdots[dotsize=0.12](9.821015,-3.0649638)
\usefont{T1}{ptm}{m}{n}
\rput(9.992471,2.7600362){$z'$}
\usefont{T1}{ptm}{m}{n}
\rput(9.822471,-3.4399638){$z$}
\usefont{T1}{ptm}{m}{n}
\rput(11.092471,2.4400363){$z'$}
\usefont{T1}{ptm}{m}{n}
\rput(12.10247,-2.2399638){$z$}
\rput{163.1795}(27.201488,-1.4916656){\psarc[linewidth=0.04,arrowsize=0.05291667cm 2.0,arrowlength=1.4,arrowinset=0.4]{->}(13.711016,1.2650362){3.11}{0.0}{78.25404}}
\usefont{T1}{ptm}{m}{n}
\rput(10.9824705,-1.6799638){$\mh^{n+1}$}
\rput{-214.79297}(24.57454,-7.1094656){\psarc[linewidth=0.04,arrowsize=0.05291667cm 2.0,arrowlength=1.4,arrowinset=0.4]{->}(13.401015,0.2950362){3.28}{0.0}{90.555504}}
\rput{-228.14397}(21.294994,-10.2231455){\psarc[linewidth=0.04,arrowsize=0.05291667cm 2.0,arrowlength=1.4,arrowinset=0.4]{->}(12.931016,-0.3549638){3.39}{0.0}{96.78173}}
\usefont{T1}{ptm}{m}{n}
\rput(11.64247,-2.6999638){$z$}
\usefont{T1}{ptm}{m}{n}
\rput(10.942471,-3.2399638){$z$}
\rput{-46.449833}(3.703228,4.8999844){\psarc[linewidth=0.04,arrowsize=0.05291667cm 2.0,arrowlength=1.4,arrowinset=0.4]{->}(7.5610156,-1.8649638){0.72}{0.0}{166.91537}}
\rput{-50.431915}(3.5755043,5.022916){\psarc[linewidth=0.04,arrowsize=0.05291667cm 2.0,arrowlength=1.4,arrowinset=0.4]{->}(7.1210155,-1.2849638){1.48}{0.0}{141.8291}}
\usefont{T1}{ptm}{m}{n}
\rput(6.8224707,0.22003618){$z$}
\usefont{T1}{ptm}{m}{n}
\rput(6.8424706,-1.2399638){$z$}
\psline[linewidth=0.04cm](4.841016,1.4950362)(3.5610156,0.37503618)
\psbezier[linewidth=0.04](3.5810156,0.36894923)(3.0410156,-0.24496381)(3.2210157,-0.60496384)(3.3410156,-1.3249638)
\psline[linewidth=0.04cm,arrowsize=0.05291667cm 2.0,arrowlength=1.4,arrowinset=0.4]{->}(3.3610156,-1.3249638)(3.9410157,-2.6849637)
\psline[linewidth=0.04cm](1.7010156,2.295036)(2.3210156,0.95503616)
\psbezier[linewidth=0.04](2.330196,0.9150362)(2.6210155,0.2950362)(2.1610155,-0.024963815)(1.6810156,-0.5649638)
\psline[linewidth=0.04cm,arrowsize=0.05291667cm 2.0,arrowlength=1.4,arrowinset=0.4]{->}(1.6610156,-0.5649638)(0.52101564,-1.7849638)
\usefont{T1}{ptm}{m}{n}
\rput(0.2224707,-1.9399638){$z$}
\usefont{T1}{ptm}{m}{n}
\rput(3.1024706,1.5600362){$g_{ij}=\del_{ij}$}
\end{pspicture} 
}
\caption{The scattering relation for compactly supported perturbations of the Euclidean metric and for hyperbolic space.}
\label{fig6N}
\end{figure}

There is a long  list of papers dedicated to the scattering relation in different settings, and it would be impossible to give a precise history.  Guillemin \cite{victor}  studied the scattering relation and sojourn times for scattering by a convex obstacle and  scattering for the automorphic wave equation after Faddeev-Lax-Phillips.  Uhlmann \cite{gunther1} showed that the Dirichlet-to-Neumann map for the wave equation gives the scattering relation on a manifold with boundary without any assumptions on caustics, while a similar result had been proved by Sylvester and Uhlmann \cite{gunther2} when there are no conjugate points.  Melrose, S\'a Barreto and Vasy \cite{MSV} studied the scattering relation for certain perturbations of the hyperbolic space. The scattering relation has also been studied in locally symmetric spaces, see for example the work of  Ji and Zworski \cite{jizw} and references cited there.   Uhlmann, Pestov and Uhlmann, Stefanov and Uhlmann \cite{pestov-uhlmann,steuhl1,steuhl2,steuhl3} studied the lens rigidity and boundary rigidity inverse problems, where one wants to obtain information about the manifold from the scattering relation.

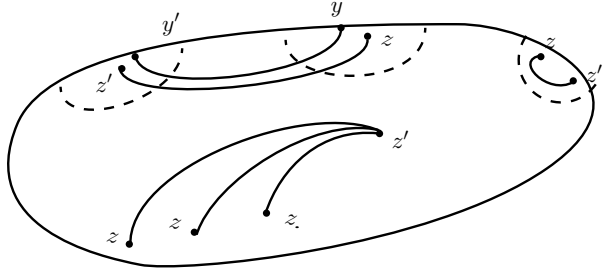
\begin{figure}
\scalebox{.8} 
{
\begin{pspicture}(0,-2.29)(11.475897,2.31)
\psdots[dotsize=0.12](6.735897,0.05)
\usefont{T1}{ptm}{m}{n}
\rput(7.0858974,-0.045){$z'$}
\psarc[linewidth=0.04](5.365897,-1.52){0.01}{0.0}{180.0}
\usefont{T1}{ptm}{m}{n}
\rput{-13.890143}(0.48434827,1.2238214){\rput(5.2655954,-1.3811944){$z$}}
\psbezier[linewidth=0.04](6.715897,0.11)(5.655897,0.35)(3.9758973,-0.73)(3.6958973,-1.63)
\psdots[dotsize=0.12](3.6558974,-1.59)
\psbezier[linewidth=0.04](6.715897,0.054285713)(5.775691,0.15)(5.093189,-0.5)(4.8358974,-1.29)
\psdots[dotsize=0.12](4.8558974,-1.27)
\usefont{T1}{ptm}{m}{n}
\rput{-13.890143}(0.44260693,0.7611497){\rput(3.3455956,-1.4411944){$z$}}
\psbezier[linewidth=0.04](8.035897,1.87)(6.7757044,1.8693372)(1.4317945,1.9271458)(0.71589726,0.21)(0.0,-1.5071459)(2.0358973,-1.99)(2.7558973,-2.13)(3.4758973,-2.27)(7.575897,-1.87)(9.515897,-0.53)(11.455897,0.81)(9.29609,1.8706628)(8.035897,1.87)
\psdots[dotsize=0.12](2.6758974,1.33)
\psdots[dotsize=0.12](6.095897,1.81)
\psbezier[linewidth=0.04](2.6758974,1.37)(2.7158973,0.59)(6.115897,1.05)(6.095897,1.79)
\psbezier[linewidth=0.04](2.4158974,1.17)(2.3958972,0.41)(6.575897,0.91)(6.5358973,1.61)
\psdots[dotsize=0.12](2.4558973,1.13)
\psdots[dotsize=0.12](6.5358973,1.67)
\usefont{T1}{ptm}{m}{n}
\rput{-13.890143}(-0.1842461,1.6950608){\rput(6.865596,1.5988057){$z$}}
\usefont{T1}{ptm}{m}{n}
\rput(2.1658974,0.895){$z'$}
\usefont{T1}{ptm}{m}{n}
\rput(3.2758973,1.875){$y'$}
\usefont{T1}{ptm}{m}{n}
\rput(6.025897,2.115){$y$}
\psbezier[linewidth=0.04](2.5958972,-1.83)(2.5158973,-0.51)(5.555897,0.63)(6.755897,0.09)
\usefont{T1}{ptm}{m}{n}
\rput{-13.890143}(0.46980968,0.5044681){\rput(2.3055956,-1.6811944){$z$}}
\psdots[dotsize=0.12](2.5758972,-1.79)
\psbezier[linewidth=0.04,linestyle=dashed,dash=0.16cm 0.16cm](1.4358972,0.83)(1.4358972,0.03)(3.4558973,0.6044262)(3.4558973,1.47)
\psbezier[linewidth=0.04,linestyle=dashed,dash=0.16cm 0.16cm](5.175897,1.61)(5.175897,0.81)(7.4958973,0.75)(7.4958973,1.81)
\rput{145.09332}(18.477665,-3.2632544){\psarc[linewidth=0.04,linestyle=dashed,dash=0.16cm 0.16cm](9.751822,1.2730958){0.7035353}{0.0}{180.0}}
\usefont{T1}{ptm}{m}{n}
\rput(9.575897,1.555){$z$}
\usefont{T1}{ptm}{m}{n}
\rput(10.305898,1.055){$z'$}
\psbezier[linewidth=0.04](9.427339,1.35)(8.935897,1.2802094)(9.53718,0.67)(9.955897,0.9089131)
\psdots[dotsize=0.12](9.415897,1.33)
\psdots[dotsize=0.12](9.955897,0.93)
\end{pspicture} 
}
\caption{Different scenarios of points approaching $\p X$: either $z'$ is fixed and $z\rightarrow \p X$ or $z, z'\rightarrow \p X,$ but $z$ and $z'$ are far apart, or $z, z'\rightarrow \p X$ and the points are close.}
\label{fig01}
\end{figure}

Let $z,z'\in \intx,$ and let $\gamma$ be a geodesic joining $z$ and $z',$ and  let $z\rightarrow y\in \p X.$  S\'a Barreto and Wunsch proved in \cite{SW} that if $\rho$ is a defining function of $\p X,$ and $\intx$ is non-trapping,  the limit $s_\gamma(z',y)=\lim_{t\rightarrow \infty} ( t+\log \rho(\gamma(t)))$ exists and moreover $s_\gamma(z',y)$ is a $C^\infty$ function of $z',y$ and the co-vector $\zeta\in T_z^*\intx$ that defines the geodesic $\gamma.$   Of course,  $s_\gamma(z',y)$ depends on the choice of $\rho.$  One of our goals here is to generalize this result to the case where both points $z$ and $z'$ joined by a geodesic $\gamma$ are allowed to go to $\p X,$ see Fig.\ref{fig01}. 

\begin{figure}[htbp]
\centering
\scalebox{0.7}
{
\begin{pspicture}(0,-3.8529167)(17.765833,3.8729167)
\rput(3.7458334,-0.14708334){\psaxes[linewidth=0.04,labels=none,ticks=none,ticksize=0.10583333cm](0,0)(0,0)(4,4)}
\rput(0.74583334,-3.1470833){\psaxes[linewidth=0.04,labels=none,ticks=none,ticksize=0.10583333cm](0,0)(0,0)(4,4)}
\psline[linewidth=0.04cm,arrowsize=0.05291667cm 2.0,arrowlength=1.4,arrowinset=0.4]{->}(3.7458334,-0.14708334)(0.24583334,-3.6470833)
\psline[linewidth=0.04cm,dotsize=0.07055555cm 3.0]{*-}(2.2458334,-1.6470833)(4.7458334,3.1729167)
\psline[linewidth=0.04cm](0.74583334,0.85291666)(3.7458334,3.8529167)
\psline[linewidth=0.04cm](4.7458334,-3.1470833)(7.7458334,-0.14708334)
\rput(2.2458334,-1.6470833){\psaxes[linewidth=0.04,arrowsize=0.05291667cm 2.0,arrowlength=1.4,arrowinset=0.4,labels=none,ticks=none,ticksize=0.10583333cm]{->}(0,0)(0,0)(5,5)}
\usefont{T1}{ptm}{m}{n}
\rput(5.295833,2.6979167){$\odiag$}
\usefont{T1}{ptm}{m}{n}
\rput(2.6,-2){$\p\odiag$}
\usefont{T1}{ptm}{m}{n}
\rput(7.4058332,-1.8420833){$x'$}
\usefont{T1}{ptm}{m}{n}
\rput(2.0358334,3.3779166){$x$}
\usefont{T1}{ptm}{m}{n}
\rput(1.35,0.5179167){$X$}
\usefont{T1}{ptm}{m}{n}
\rput(4.1658335,-2.6820834){$X$}
\usefont{T1}{ptm}{m}{n}
\rput(0.97583336,-3.6020834){$y - y'$}
\rput(13.745833,-0.14708334){\psaxes[linewidth=0.04,labels=none,ticks=none](0,0)(0,0)(4,4)}
\rput(10.745833,-3.1470833){\psaxes[linewidth=0.04,labels=none,ticks=none](0,0)(0,0)(4,4)}
\psline[linewidth=0.04cm](13.745833,3.8529167)(10.745833,0.85291666)
\psline[linewidth=0.04cm](17.745832,-0.14708334)(14.745833,-3.1470833)
\psline[linewidth=0.04cm](11.865833,-2.0270834)(10.745833,-3.1470833)
\psarc[linewidth=0.04](13.425834,-2.0270834){1.56}{90.0}{180.0}
\psarc[linewidth=0.04](11.885834,-0.50708336){1.54}{268.5312}{1.27303}
\psline[linewidth=0.04cm](13.745833,-0.14708334)(13.425834,-0.48708335)
\psarc[linewidth=0.04](12.155833,-1.9770833){1.03}{28.855661}{82.11686}
\psline[linewidth=0.04cm,dotsize=0.07055555cm 3.0]{*-}(12.765833,-1.1470833)(14.725833,3.1129167)
\usefont{T1}{ptm}{m}{n}
\rput(15.245833,2.6579165){$\diag_0$}
\usefont{T1}{ptm}{m}{n}
\rput(14.445833,-1.9970833){\large $L$}
\usefont{T1}{ptm}{m}{n}
\rput(11.635834,0.8229167){\large $R$}
\usefont{T1}{ptm}{m}{n}
\rput(12.325833,-1.6020833){$\ff$}
\psline[linewidth=0.04cm,arrowsize=0.05291667cm 2.26,arrowlength=1.4,arrowinset=0.4]{->}(10.1258335,3.4129167)(8.105833,3.4129167)
\usefont{T1}{ptm}{m}{n}
\rput(9.155833,2.9029167){\large $\beta_0$}
\end{pspicture}
}
\caption{The $0$-blown-up space $X\times_0X$.}
\label{fig1}
\end{figure}
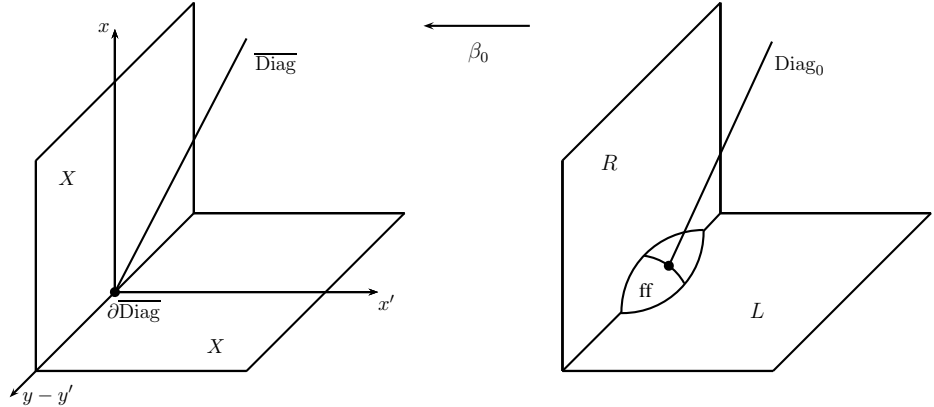

The metric $g$ on $T \intx$  induces a dual metric $g^*$ on $T^*\intx$ defined in \eqref{dual-met}.  We shall  view a geodesic as the projection of an integral curve of $H_p,$ the Hamilton vector field of $p(z,\zeta)=\ha(|\zeta|_{g^*(z)}-1),$ see Section \ref{sojourn}.  The integral curve of $H_p$ connects two points $(z',\zeta')$ and $(z,\zeta),$ and to better understand the map $(z',\zeta')\mapsto (z,\zeta),$ it will be convenient to work on the product $\intx\times \intx.$  
  We can identify $T^*(\intx \times \intx)= T^* \intx \times T^*\intx$ and  according to this, we shall  use  $(z,\zeta,z',\zeta')$ to denote a point in $T^*(\intx \times \intx),$ while $(z,\zeta)$ will denote a point on the left factor and $(z',\zeta')$ will denote a point on the right factor.  In fact we shall work on 
  $T^*(\mr\times \intx \times \intx),$  and we denote  
  \begin{gather*}
 Q_R(\tau, z,\zeta,z',\zeta')=\ha(|\zeta'|_{g^*(z')}-\tau^2) \text{  and  } Q_L(\tau,z,\zeta,z',\zeta')=\ha(|\zeta|_{g^*(z)}-\tau^2).
 \end{gather*}
    In Section \ref{sojourn}, we study the Lagrangian submanifold $\tilde \La \subset T^*(\mbr_t\times \intx\times \intx)$ defined as the flow out of 
 \begin{gather}
 \{(0,\tau,z,\zeta, z',\zeta'):  z=z', \zeta=-\zeta', \; \tau^2=|\zeta|_{g^*(z)} \}, \label{diagonal}
\end{gather}
under the Hamilton vector field $H_{Q_R}$ or $H_{Q_L}.$  Since $\intx$ is non-trapping, $\tilde\La$ is a $C^\infty$ Lagrangian submanifold of $T^*(\mr_t \times \intx \times \intx),$ and we will analyze the global behavior of $\tilde\La$ up to $\p(\mbr_t \times X \times  X).$    The obvious problem is when  the closure of $\diag=\{(z,z'):z=z', z\in \intx \}$  meets $\p X \times \p X,$ where the Hamilton flow is not well-defined.    To handle this situation, we  work in  the $0$-blow-up of $X\times X$ defined in Mazzeo and Melrose \cite{MM}, and we recall their construction.   Let 
\beq
\p \odiag = \{(z, z) \in \p X\times \p X\} = \odiag \cap(\p X\times \p X).
\eeq
As a set, the $0$-blown-up space is 
\beq
X\times_0 X = (X\times X)\backslash \p\odiag \sqcup S_{++}(\p \odiag),
\eeq
where $S_{++}(\p\odiag)$ denotes the inward pointing spherical bundle of $T_{\p\odiag}^*(X\times X)$. Let 
\beq
\beta_0: X\times_0 X \rightarrow X\times X
\eeq 
be the blow-down map. Then $X\times_0 X$ is equipped with the topology and smooth structure  of a manifold with corners such that $\beta_0$ is smooth. The manifold $\xo$  has three boundary hypersurfaces: the left and right faces $L=\overline{\beta_0^{-1}(\p X\times \intx)},$  $R=\overline{\beta_0^{-1}(\intx \times \p X)},$  and the front face $\ff= \overline{\beta_0^{-1}(\p\odiag)}$. The lifted diagonal is denoted by $\diag_0=\overline{\beta_0^{-1}(\diag)}$. See Figure \ref{fig1}.  

In the interior of $\xo,$ $\beta_0$ is a diffeomorphism between open $C^\infty$ manifolds, and $\beta_0^*{\tilde\La}$ is naturally well-defined as the joint flow-out of the lift of \eqref{diagonal} under the lifts $\beta_0^*H_{Q_R}$ and $\beta_0^* H_{Q_L}.$     By abuse of notation, we will also denote

\begin{gather*}
\beta_0: \mr_t \times \xo \longrightarrow \mr_t\times X \times X \\
(t, m) \longmapsto (t, \beta_0(m)).
\end{gather*}

 We will prove the following result in Section \ref{sojourn}:
\begin{theorem}\label{soj}
Let $(\intx, g)$ be a non-trapping AHM.  Let $\rho_L, \rho_R$ be  boundary defining functions of $L$ and $R$ respectively.  Let 
\begin{gather}
\begin{gathered}
\mcm: \mr\times \left(\xo \setminus (R\cup L)\right) \longrightarrow \mr \times \xo \\
(t, m) \longmapsto ( t + \log \rho_L(m) + \log \rho_R(m), m)=(s,m)
\end{gathered}\label{sch}
\end{gather}
and define $\beta_1=\beta_0\circ \mcm: \mr_s \times \left(\xo \setminus(R\cup L)\right) \longrightarrow \mr_t\times X \times X.$  Let  $\tilde \La \subset T^*(\mbr_t\times \intx\times\intx)$ be the $C^\infty$ Lagrangian submanifold defined in \eqref{eqlat}. Let $\beta_1^*\tilde\La$ denote the lift of $\tilde \La$ by $\beta_1$ in the interior of $\xo.$
 Then  $\beta_1^*\tilde \La$ has a smooth extension up to the boundary of  $T^*(\mbr_s \times \xo).$
\end{theorem}

The main point  in the proof of Theorem \ref{soj} is that $\diag_0$ does not intersect $R$ or $L$ and intersects $\ff$ transversally, see Fig. \ref{fig1}.  We will also show that after the singular change of variables \eqref{sch},  $H_{\frac{1}{\rho_R}\beta_1^*Q_R}$ and $H_{\frac{1}{\rho_L}\beta_1^*Q_L}$ lift to  $C^\infty$ vector fields  which are tangent to $\mr_s \times \ff.$ Morever  if $\sigma$ is the dual variable to $s,$ in the region $\sigma\not=0,$ where $\beta_1^*\tilde\La$ is contained, $H_{\frac{1}{\rho_R}\beta_1^*Q_R}$ is  transversal to  $\mr_s\times R,$
and  $H_{\frac{1}{\rho_L}\beta_1^*Q_L}$ is  transversal to  $\mr_s\times L.$  Therefore, the manifold $\beta^*\tilde\La$ extends smoothy up to 
$\p(\mr_s \times \xo),$ see Fig. \ref{fig4}.  The main point is that after the singular change of variables and rescaling  the lift of the symbols, the Lagrangian results from the integration of a $C^\infty$ vector field over a finite  interval.  This idea is reminiscent from the work of S\'a Barreto and Wunsch \cite{SW} and  Melrose, S\'a Barreto and Vasy \cite{MSV}.
 
 Theorem \ref{soj} generalizes a result of \cite{SW} and shows that one can define the sojourn time along a geodesic joining two points $z$ and $z',$ as both points go to $\p X.$ This will be discussed in details in Section \ref{sojourn}.

 As an application of Theorem \ref{soj},  we study the asymptotics of the distance function $r(z,z')$ between $z,z'\in \intx$ as $z,z'\rightarrow \p X,$ in the case where $(\intx,g)$ is a  geodesically convex AHM. In this case there are no conjugate points  along any geodesic in $\intx$ and  $r(z,z')$ is equal to the length of the unique geodesic joining the two points. Moreover, $r(z,z')$  is smooth on $(\intx \times \intx)\backslash\diag$.  This is the case when $(\intx,g)$ is a Cartan-Hadamard manifold, i.e.
 when $\intx$ has non-positive sectional curvature, see \cite{GHL}.     
 We will show in  Section \ref{distfunc} that the following is a consequence of Theorem \ref{soj}
\begin{theorem}\label{main}
Let $(\intx, g)$ be a geodesically convex AHM, and let $\rho_L$ and $\rho_R$ be boundary defining functions of  $L$ and $R$  respectively.  For $z,z' \in \intx$,  the lift of the distance function $r(z, z')$  to $\xo$ satisfies 
\begin{gather}
\beta_0^*r = -\log \rho_L-\log \rho_R + F, \phan F\in C^{\infty}(\xo\backslash \diag_0). \label{eqfr}
\end{gather}
\end{theorem}

One should emphasize that the importance of \eqref{eqfr} is what it reveals about the behavior of $\beta_0^*r$ near the right, left and front faces, and where $\diag_0$ meets $\ff.$  The singularity of $\beta_0^*r$ at the diagonal in the interior is well understood. In fact, since $r(z,z')$ is the distance function, then in the interior of $\xo$ and near $\diag_0,$ $(\beta_0^*r(z,z'))^2$ is $C^\infty$ and vanishes quadratically at $\diag_0.$    One can always modify $\rho_R$ and $\rho_L$  such that $\rho_R=\rho_L=1$ near $\diag_0,$ and  with this choice of $\rho_R$ and $\rho_L$ , one has that $F=\beta_0^*r$ near $\diag_0.$  This regularity near $\diag_0$ extends up to $\ff.$

The asymptotics of $r(z,z')$  for perturbations of the hyperbolic space  of the type
\begin{gather}
\intx=\mb^{n+1} \text{ equipped with the metric } g_\eps=\frac{4 dz^2}{(1-|z|^2)^2} + \chi(\frac{1-|z|^2}{\eps}) H(z,dz), \label{msv-model}
\end{gather}
where $\chi(t)\in C_0^\infty(\mbr),$ with $\chi(t)=1$ if $|t|<1$ and  $\chi(t)=0$ if $|t|>2,$  $H$ is a $C^\infty$ symmetric 2-tensor and $\eps$ is small enough, was studied by Melrose, S\'a Barreto and Vasy in \cite{MSV}  in connection with the analysis of the asymptotic behavior of solutions of the wave equation on de Sitter-Schwarszchild space-time.    For $\eps$ small enough, $(\mb^{n+1}, g_\eps)$ is an example of a Cartan-Hadamard manifold.  It was proved in \cite{MSV}  that there exists $\eps_0>0$ such that if $\eps\in (0,\eps_0),$ then Theorem \ref{main} holds for the particular case when $g$ is given by \eqref{msv-model}.  This was proved by first  analyzing the case of the hyperbolic space and using  perturbation arguments. 

Equation \eqref{eqfr} was the key ingredient in the construction of  a semiclassical parametrix for the resolvent of the Laplacian
 $R(\la,h)=\left(h^2(\Delta_g-\frac{n^2}{4})-\la^2\right)^{-1},$ when $(\intx,g)$ satisfied \eqref{msv-model} carried out in \cite{MSV}.     In particular, in view of \eqref{eqfr}, the proof of Theorem 5.1 of \cite{MSV} extends line by line to geodesically convex  AHM.

 It is easy to illustrate Theorem \ref{main}  in half-space model of  hyperbolic space. A similar computation is done in \cite{MSV}.  In this case,  
\begin{gather*}
\mh^{n+1}=\{(x,y): \ \ x>0, \;\ y\in \mr^n\}, \;\ g=\frac{dx^2}{x^2}+ \frac{dy^2}{x^2},
\end{gather*}
and the distance function satisfies
\begin{gather*}
\cosh r(z,z')=\frac{x^2+ {x'}^2+|y-y'|^2}{2x x'}, \;\ z=(x,y), \;\ z'=(x',y').
\end{gather*}

Since
\begin{equation}\label{eqer}
e^r = \cosh r + \sqrt{\cosh^2 r - 1},
\end{equation}
we obtain
\begin{gather*}
{r(z,z')}= -\log(x x') +\\ \log\ha\left( x^2+x'^2+ |y-y'|^2 + \left[ \left((x+x')^2+|y-y'|^2\right)\left((x-x')^2+|y-y'|^2\right)\right]^\ha\right).
\end{gather*}

This shows that away from the diagonal, i.e. if  $(x-x')^2+|y-y'|^2>\del>0,$
\begin{gather*}
r(z,z')=-\log(x x') + F,
\end{gather*}
where $F$ is smooth up to $x = 0$ and $x' = 0$. 
However the asymptotic behavior of  $r(z,z')$ near  the diagonal as $x\downarrow 0$ and $x'\downarrow 0$  is  more appropriately expressed in terms of polar coordinates.  In this case, we choose $R=[x^2+{x'}^2+|y-y'|^2]^\ha$ as a defining function of the submanifold
\beq 
\p \odiag= \{x=0, \; x'=0, \ \ y = y'\},
\eeq
 and  denote $\rho_R=x'/R,$ $\rho_L=x/R$ and $Y=\frac{y-y'}{R}.$  Then by using \eqref{eqer}, we have
\begin{gather}
\beta_0^* r +\log (\rho_R\rho_L) = \log \ha\left(1+ \left[ \left((\rho_R+\rho_L)^2+|Y|^2\right)\left( (\rho_R-\rho_L)^2+|Y|^2\right)\right]^\ha\right),\label{asymodel1}
\end{gather}
which confirms \eqref{eqfr} in this example.

Notice that the right hand side of \eqref{asymodel1} does not depend on $R.$ So the asymptotics holds up to the front face $\ff.$ Notice also that, away from the front face
\begin{gather*}
r(z,z')+\log(x x')|_{\{x=x'=0\}}= 2\log|y-y'|.
\end{gather*}
This gives the sojourn time between two points that approach $\p X$ away from each other, see Fig.\ref{fig01}. Similarly,
\begin{gather}
\beta_0^*( r+ \log x + \log x')|_{\{\rho_R=\rho_L=0\}}=2\log(R|Y|)|_{\{\rho_R=\rho_L=0\}}= 2\log R|_{\{\rho_R=\rho_L=0\}}, \label{SOJT}
\end{gather}
since $|Y|=1$ when $\rho_L=\rho_R=0.$     It is important to point out the difference between the two asymptotic expansions. From \eqref{asymodel1} we have
\begin{gather*}
(\beta_0^* r+ \log \rho_R+ \log \rho_L)|_{\{\rho_R=\rho_L=0\}}=0.
\end{gather*}
According to \cite{SW}, the sojourn time between points on $\xo$ is given by \eqref{SOJT}.  Notice that, according to \cite{JS1} this is precisely the singularity of the Schwartz kernel of the scattering matrix.

\section{The Lagrangian manifold, the scattering relation and sojourn times}\label{sojourn}
We recall some basic facts about Riemannian and symplectic geometry. The Riemannian metric $g$ in the interior $\intx$ induces an isomorphism
\begin{gather}
\begin{gathered}
\mathcal{G}: T_z\intx\longrightarrow T_z^* \intx \\
v\longmapsto g(z)(v, \cdot),
\end{gathered}\label{isomor}
\end{gather}
which in turn induces a dual metric on $T^* \intx$ given by
\begin{gather}
g(z)^*(\xi,\eta)= g(z)( \mcg^{-1} \xi, \mcg^{-1}\eta). \label{dual-met}
\end{gather}
In local coordinates we have
\begin{gather*}
g(z)(v,w)=\sum_{i,j} g_{ij}(z) v_i w_j \text{ and }  g(z)^*(\xi,\eta)=\sum_{i,j} g^{ij} (z)\xi_i \eta_j, \text{ where  the matrices } (g_{ij})^{-1}=(g^{ij}).
\end{gather*}

Consider the product manifold $T^*(\mr_t \times \intx\times\intx)$ which can be identified with $T^*\mr \times T^*\intx\times T^*\intx$. 
In local coordinates $(z, \zeta, z', \zeta')$  on $T^*(\intx\times\intx)$, the canonical $2$-form is
\beq
\tilde\omega = d\tau \wedge dt+ \sum_{j = 1}^{n+1} d\zeta_j\wedge dz_j + \sum_{j = 1}^{n+1} d\zeta'_j\wedge dz'_j.
\eeq
 We shall denote the conormal bundle of the diagonal by 
\begin{gather*}
N^* \diag\setminus 0=\{(z,\zeta,z',\zeta'): z'=z\in \intx, \;\ \zeta=-\zeta'\not=0\}.
\end{gather*}
We distinguish between the lift of the wave equation associated to the right or left factor.  Let
\begin{gather*}
\square_{g_R}= \ha(D_t^2-\Delta_{g_R}) \text{ and }
\square_{g_L}= \ha(D_t^2-\Delta_{g_L}).
\end{gather*}
The principal symbols of these operators are 
\begin{gather}
\begin{gathered}
 Q_R(z', \zeta', t, \tau)=\ha( \tau^2- |\zeta'|_{g^*(z')}^2), \text{ and } \\
 Q_L(z, \zeta, t, \tau)=\ha( \tau^2- |\zeta|_{g^*(z)}^2).
\end{gathered}\label{eqqr}
\end{gather}
 In what follows we will think of these as functions on $T^*(\mbr\times\intx\times \intx),$ and $H_{Q_R}$ and $H_{Q_L}$ will denote
  the Hamilton vector fields of $Q_R$ and $Q_L$ with respect to $\tilde{\omega}.$   In local coordinates
  \begin{gather*}
  H_{Q_\bullet}= \tau\frac{\p}{\p t} +  \sum_{j=1}^{n+1} \frac{\p Q_\bullet}{\p {\zeta_j}}\frac{\p}{\p{z_j}}-\frac{\p Q_\bullet}{ \p{z_j}}\frac{\p}{ \p{\zeta_j}}, \text{ where } \bullet=R, L.
  \end{gather*}
 These vector fields obviously commute, and hence for $t_1\geq 0$ and $t_2 \geq 0,$ and $(t,\tau,z,\zeta,z',\zeta')\in T^*(\mbr\times\intx\times\intx),$
\begin{gather*}
\exp t_2 H_{Q_L}\circ \exp t_1 H_{Q_R}(t,\tau,z,\zeta,z',\zeta')= \exp t_2 H_{Q_R}\circ \exp t_1 H_{Q_L}(t,\tau,z,\zeta,z',\zeta').
\end{gather*}

With this identification, if
\begin{equation}\label{initial}
\Sigma \doteq  \{ t=0, \; z=z'  \in \intx, \;  \zeta=-\zeta'\not=0, \;\ \tau^2=|\zeta|_{g^*(z)}^2\},
\end{equation}
we define
\begin{gather}\label{eqlat}
\tilde\La \doteq \bigcup_{t\geq 0} \exp t H_{Q_R} (\Sigma)= \bigcup_{t \geq 0} \exp t H_{Q_L} (\Sigma)= \bigcup_{t_1\geq 0, t_2\geq 0} \exp t_2 H_{Q_L}\circ \exp t_1 H_{Q_R} (\Sigma).
\end{gather}
To see the last equality, we  just have to realize that  if $(t,\tau,z,\zeta)=\exp(sH_{Q_R})(t_1,\tau,z',\zeta'),$ then
$(t_1,\tau,z',\zeta')=\exp(s H_{Q_L})(t,\tau,z,\zeta).$

Now we will carry out the proof of Theorem \ref{soj}. First, notice that the result is independent of the choice of $\rho_R$ or $\rho_L.$   If  $\tilde\rho_L, \tilde \rho_R$ are  boundary defining functions of the left and right faces, then $\rho_L = \tilde \rho_L f_L$ and $\rho_R = \tilde \rho_R f_R$ for some $f_L, f_R \in C^\infty(\xo)$ with $f_L>0,$ $f_R > 0.$ If  $\tilde{s}= t + \log \tilde \rho_L + \log \tilde \rho_R$ and $s=t + \log\rho_R+\log\rho_L,$ then $\tilde{s}=  s + \log(f_L f_R),$ and the map $(s,m)\mapsto(\tilde{s},m)$ is a global diffeomorphism of $\mr_s\times \xo.$

 The main ingredient in the proof of Theorem \ref{soj} is the following
 \begin{lemma}  Let $\rho_R, \rho_L\in C^\infty(\xo)$ be defining functions of  $R$ and $L$ respectively.   Let $\beta_1$ be the map defined in \eqref{sch}
  and let   $q_R=-\frac{1}{\rho_R}\beta_1^*Q_R$ and
$q_L= -\frac{1}{\rho_L} \beta_1^* Q_L.$ Then $q_R$ and $q_L$ extend to  functions in $C^\infty(T^*(\mr_s\times \xo))$ and the Hamilton vector fields  $H_{q_L}$ and $H_{q_R}$ are tangent to $\mr_s\times \ff.$ Moreover, if $\sigma$ is the dual variable to $s,$ then away from $\sigma=0,$ $H_{q_R}$ is transversal to $\mr_s\times R,$ and $H_{q_L}$ is transversal to $\mr_s\times L.$
 \end{lemma}
 \bpf   We will prove this Lemma in local coordinates valid near $\p(\mr_s \times \xo).$   First, we choose local coordinates  $z = (x, y)$ and $z' = (x', y')$ in which \eqref{prod} holds.  Then we pick the following defining functions of $\ff,$ $R$ and $L:$ 
\begin{gather}
\begin{gathered}
\rho_{\ff}=\left[ x^2+(x')^2+|y-y'|^2\right]^\ha \text{ is a defining function of } \ff \\
\rho_R=\frac{x'}{\rho_{\ff}} \text{ is a defining function of  the right face } R \\
\rho_L=\frac{x}{\rho_{\ff}} \text{ is a defining function of  the left face }  L.
\end{gathered}\label{coords}
\end{gather}

We will divide the boundary of $\mr_s\times \xo$ in four regions: \\
Region 1:  Near $\mr_s \times L,$ and away from $\mr_s\times (R\cup \ff),$ or near $\mr_s\times R$ and away from $\mr_s\times (L \cup \ff).$ \\
Region 2: Near $\mbr_s\times (L \cap \ff)$ and away from $\mbr_s\times R$, or near $\mbr_s\times (R\cap \ff)$ and away from $\mbr_s\times L.$ \\
Region 3:  Near $\mbr_s\times (L\cap R)$ but away from $\mbr_s\times \ff.$ \\
Region 4: Near $\mbr_s\times (L\cap R \cap \ff).$

  First we analyze region 1, near $\mr_s\times L $ but away from $\mr_s \times R$ and $\mr_s \times \ff.$  The case near $\mr_s\times R$ but away from $\mr_s\times L$ and $\mr_s\times \ff$ is identical.   Since we are away from $R,$ we have $\rho_R>\del,$ for some $\del>0,$  and hence $\log \rho_R$ is $C^\infty.$ In this region we may take $x$ as a defining function of $L,$ and instead of \eqref{sch}, we set $s = t + \log x$.   In fact,   the map $(s,m) \longmapsto (s+\log\rho_R,m)$ is a diffeomorphism in the region where $\rho_R>\del,$ and hence  the statements about $q_L$ and $H_{q_L}$ in the lemma are true in this region whether we take $s=t+\log x$ or $s=t+\log x + \log \rho_R.$  In the case near $\mr_s\times R$ but away from $\mr_s\times L$ and $\mr_s\times  \ff$ one sets $s=t+\log x'.$  These particular cases were studied in  \cite{SW}.
  
The change of variables
\begin{gather}
s=t+\log x \label{souj1}
\end{gather}
induces the symplectic change on $T^*(\mbr\times \intx\times \intx)$

\begin{gather}\label{symch}
\begin{gathered}
(x,y,\xi,\eta,t,\tau) \longmapsto (x,y,\txi,\eta,s,\sigma), \\
\text{ where } \txi= \xi-\frac{1}{x}\tau, \ \ \sigma = \tau.
 \end{gathered}
  \end{gather}
  
 In coordinates \eqref{prod}, 
 \begin{gather*}
 g_L^*(x,y,\xi,\eta)= x^2 \xi^2+ x^2 h(x,y,\eta),
 \end{gather*}
 and so
 \begin{gather*}
 \beta_1^*Q_L=-x\sigma\txi-\ha x^2(\txi^2+h(x,y,\eta)), \text{ and hence } \\
   q_L= -\frac{1}{\rho_L} \beta_1^*Q_L=\sigma\txi+\ha x(\txi^2+h(x,y,\eta)).
 \end{gather*}

We have
 \begin{gather*}
 H_{q_L}= (\sigma+x\txi)\p_x+\txi\p_s +\ha x H_{h(x,y,\eta)}- \ha(\txi^2+h(x,y,\eta)+ x\p_x h(x,y,\eta)) \p_{\txi}.
 \end{gather*}
 In particular, $\sigma$ remains constant along the integral curves of $H_{q_R},$ and
 \begin{gather*}
H_{q_L}|_{\{x=0\}}= \sigma \p_x+\txi\p_s - \ha(\txi^2 + h(0, y, \eta))\p_{\txi}.
\end{gather*}
So if $\sigma\not=0,$  $H_{q_L}$ is transversal to $\p X.$

Next we work in region 2 near $\mr_s \times (L\cap \ff),$ but away from $\mr_s\times R$.  The case near $\mr_s\times (R\cap \ff)$ but away from $\mr_s\times L$ is very similar.   In this case, $\rho_R=x'/R>\del,$ and so it is better to use projective coordinates
\begin{equation}\label{eqc1}
X = \frac{x}{x'},\ \ Y = \frac{y - y'}{x'},\ \ x' \text{ and } y'. 
\end{equation}
In this case, $X$ is a boundary defining function for $L$ and $x'$ is a boundary defining function for $\ff$. Since $\beta_0$ is a diffeomorphism in the interior of $X\times_0 X$, it induces a symplectic change of variables 
\beq
(x, y, \xi,\eta, x', y', \xi', \eta')\in T^*(\intx \times \intx) \longmapsto (X, Y, \la,\mu, x', y', \la', \mu') \in T^*(X\times_0 X),
\eeq
given by
\beq
\la = x' \xi,\ \ \mu = x'\eta,\ \ \la' = \xi' + \xi X + \eta Y  \text{ and } \mu' = \eta + \eta'.
\eeq
and  $Q_L$  becomes
\beq
\beta_0^*Q_{L}  = \ha(\tau^2 - X^2( \la^2 + h(x'X, x'Y + y', \mu))),
\eeq
and here we used the fact that $h(x,y,\eta)$ is homogeneous of degree two in $\eta.$

  Away from the face $R$, $\rho_R>\del,$ for some $\del,$ and the function $\log \rho_R$ is smooth. Therefore, as argued above in the case of region 1, the transformation  $(s,m) \mapsto (s+\log \rho_R,m)$ is a $C^\infty$ map away from $\{\rho_R=0\},$ and so it suffices to take
\begin{equation}\label{ch1}
s = t + \log X.
\end{equation}

The change of variable \eqref{ch1}  induces the following symplectic change of variables
\begin{gather}
\begin{gathered}
T^*(\mbr_t \times \intx \times \intx) \longrightarrow  T^*(\mbr_s \times X\times_0 X),\\
(t, \tau, x, y,\xi,\eta, x', y', \xi', \eta')  \longmapsto (s,\sigma, X, Y,\tilde\lambda,\mu, x', y', \la', \mu')\\
\text{where } \tilde \la = \la - \frac{\tau}{X}, \ \ \sigma = \tau,
\end{gathered}
\end{gather}
and the canonical $2$-form on $T^*(\mbr_s\times\xo)$ is given by
\beq
\omega^0 = d\tilde\la\wedge d X + d\mu \wedge dY + d\la' \wedge dx' + d\mu'\wedge dy'.
\eeq
Hence
\beq
\beta_1^* Q_{ L} = - \tilde\la\sigma X - \ha X^2 ({\tilde\la}^2 + h(x'X, x'Y + y', \mu)),
\eeq
and we conclude that
\beq
q_{L} = -\frac{1}{\rho_L}\beta_1^* Q_L=-\frac{1}{X}\beta_1^*Q_L = \tilde\la\sigma + \ha X (\tilde\la^2 + h(x'X, x'Y + y', \mu)).
\eeq
Hence vector field $H_{q_{L}}$ is given by
\begin{gather}
\begin{gathered}
H_{q_L} = \tilde\la \frac{\p}{\p s} + (\sigma + X \tilde\la)\frac{\p }{\p X}  - \ha(\tilde\la^2 + h + x'X\p_X h)\frac{\p}{\p \tilde\la} + \frac{X}{2} H_h + T,
\end{gathered}\label{liftedHA}
\end{gather}
where $T$ is a smooth vector field in $\p_{\la'}, \p_{\mu'}$.   So away from $\sigma=0,$ $H_{q_L}$ is transversal to $\mr_s\times L.$

Next we analyze region 3, near $\mr_s\times (L\cap R)$ and away from $\mr_s\times \ff$.  Here $x, x'$ are boundary defining functions for $\mr_s\times L$ and  $\mr_s\times R$ respectively. In this case, as discussed above, we can take
\beq
s = t + \log x + \log x',
\eeq
which induces the following symplectic change of variable
\begin{gather*}
(t, \tau,x, y,\xi,\eta, x', y' , \xi', \eta') \longmapsto (s,\sigma, x, y, \tilde\xi, \eta, x', y', \tilde\xi', \eta'),\\
\text{where } \tilde \xi =  \xi - \frac{\tau}{x}, \ \ \tilde \xi' = \xi' - \frac{\tau}{x'}, \ \ \sigma = \tau.
\end{gather*}
The symbols can be computed as in the case near $\mr_s\times L$ away from $\mr_s\times \ff$ and $\mr_s\times R$. In particular,
\begin{gather*}
\beta_1^*Q_L=-x\sigma\txi-\ha x^2(\txi^2+h(x,y,\eta))  \text{ and so } q_L=-\frac{1}{\rho_L} \beta_1^* Q_L=\sigma\txi+\ha x(\txi^2+h(x,y,\eta)),\\
\beta_1^*Q_R=-x\sigma\txi'-\ha x'^2(\txi'^2+h(x',y',\eta'))  \text{ and so }  q_R=\sigma\txi'+\ha x'(\txi'^2+h(x',y',\eta')).
\end{gather*}

The Hamilton vector fields are given by
\beq
\bsp
& H_{q_L}= (\sigma+x\txi)\p_x+\txi\p_s +\ha x H_{h(x,y,\eta)}- \ha(\txi^2+h(x,y,\eta)+ x\p_x h(x,y,\eta)) \p_{\txi}, \\
& H_{q_R}= (\sigma+x'\txi')\p_{x'}+\txi'\p_s +\ha x' H_{h(x',y',\eta')}- \ha(\txi'^2+h(x',y',\eta')+ x'\p_{x'} h(x',y',\eta')) \p_{\txi'}.
\end{split}
\eeq 
We conclude that, away from $\sigma=0,$  $H_{q_L}$ is transversal to $\mr_s\times L=\{x = 0\}$ while $H_{q_R}$ is transversal to $\mr_s\times R=\{x' = 0\}.$ 

Finally, we analyze region 4,  near  the co-dimension $3$ corner $\mr_s\times (L\cap \ff\cap R).$  Here we also work with suitable projective coordinates, and without loss of generality, as in \cite{MSV} we may take $\rho_{\ff}=y_1 - y_1' \geq 0$ and take the following coordinates
\begin{equation}\label{eqc2}
u = y_1 - y_1', \ \ w = \frac{x}{y_1 - y_1'}, \ \ w' = \frac{x'}{y_1 - y_1'},\ \ y'  \text{ and } Z_j = \frac{y_j - y_j'}{y_1 - y_1'}, \ \ j = 2, 3,\cdots n.
\end{equation}
Here $w, w'$ and $ u$ are boundary defining functions for $\mr_s\times L,$  $\mr_s\times R$ and $\mr_s\times \ff$ faces respectively.  The induced symplectic change of variables
\begin{gather}
\begin{gathered}
T^*(\intx\times\intx) \longrightarrow  T^*(X\times_0 X) \\ 
(x, y, \xi,\eta, x', y', \xi', \eta') \longmapsto (w,u,Z,\la,\nu,\mu, w', y',  \la', \mu') \\
\text{ where } \\
\la  = \xi u, \ \ \la'  = \xi' u, \ \ \nu = \xi w + \xi'w' + \eta_1 + \sum_{j=2}^n \eta_j Z_j, \\
\mu' = \eta + \eta',  \ \ \mu_j = \eta_j u, \ \ j = 2, 3, \cdots n. 
\end{gathered}\label{eqsym}
\end{gather}

In these coordinates, the symbols of $Q_L$ and $Q_R$ are given by
\begin{gather*}
\beta_0^* Q_{ L}  = \ha( \tau^2 -  w^2(\la^2 + h(uw, y, u\eta))),\\
 \beta_0^* Q_{ R}  = \ha( \tau^2 -  w'^2(\la^2 + h(uw', y', u\mu' - u\eta))),
\end{gather*}
where 
\beq
y = (y_1' + u, y_2' + uZ_2, \cdots, y_{n}' + uZ_{n}),\ \ u\eta = (u\nu - \la w - \la' w' - \sum_{j = 2}^n \mu_j Z_j, \mu).
\eeq
In this case,  we set
\begin{equation}\label{ch2}
s = t + \log w + \log w',
\end{equation}
which  induces the symplectic transformation
\beq
\bsp
& T^*(\mbr_t \times \intx\times\intx) \longrightarrow  T^*(\mbr_s \times X\times_0 X),\\
& (t, \tau,x, y, \xi,\eta, x', y',  \xi', \eta')  \longmapsto (s,\sigma, w,u,Z, \tilde\la, \nu,\mu, w', y', \tilde\la',  \mu') \\
&  \text{where } \tilde \la = \la - \frac{\tau}{w}, \ \ \tilde \la' = \la'- \frac{\tau}{w'}, \ \ \sigma = \tau.
\end{split}
\eeq
Here the canonical $2$-form on $T^*(\mbr_s \times \xo)$ is given by 
\beq
\omega^0 = d\sigma \wedge ds + d\tilde\la \wedge dw + d\tilde\la'\wedge dw' + d\nu\wedge du + d\mu\wedge dZ + d\mu' \wedge d y'.
\eeq
The lift of the symbols $Q_{L}$ and $Q_{R}$ become
\begin{gather*}
\beta_1^*Q_{ L} = -w\sigma \tilde \la - \ha w^2 (\tilde\la^2 + h(uw, y, \tilde \eta)),\\
 \beta_1^*Q_{ R} = -w'\sigma \tilde\la' - \ha w'^2 ( \tilde\la'^2 + h(uw', y', u\mu' - \tilde \eta)),
\end{gather*}
where $\tilde \eta \doteq u\eta =  (u\nu - \tilde\la w - \tilde\la' w' - 2\sigma - \sum_{j = 2}^n \mu_j Z_j, \mu)$.  Therefore, in these coordinates
\begin{gather*}
q_{L} = -\frac{1}{\rho_L}\beta_1^*Q_L =  \sigma \tilde\la + \frac{1}{2} w(\tilde\la^2 +  h(uw,  y, \tilde \eta)),\\
q_{R} =-\frac{1}{\rho_R}\beta_1^*Q_R=  \sigma \tilde\la' + \frac{1}{2}w '( \tilde\la'^2 +  h(uw', y', u\mu' - \tilde \eta)).
\end{gather*}
Hence the Hamilton vector fields are of the form
\begin{gather}
\begin{gathered}
H_{q_{ L}} = (\sigma + w\tilde\la - \ha w^2 \p_{\tilde\eta}h(uw, y, \tilde \eta))\frac{\p}{\p w} + T_L,\\
  H_{q_{R}} =  (\sigma + w'\tilde\la' + \ha w'^2 \p_{\tilde\eta}h(uw', y', u\mu'- \tilde \eta) )\frac{\p}{\p w'} + T_R,
\end{gathered}\label{VF}
\end{gather}
 where $T_L, T_R$ are smooth vector fields on $T^*(\mbr_s \times X\times_0 X)$ with no $\frac{\p}{\p w},$ $\frac{\p}{\p w'}$ or $\frac{\p}{\p \sigma}$  terms. 
Notice that these vector fields are $C^\infty$ up to the front face, and that away from $\sigma=0,$ the vector $H_{q_L}$ is transversal to $\mr_s\times L$  and $H_{q_R}$ is transversal to $\mr_s\times R.$  This shows that the transversality to $L$ and $R$ holds up to the corner.
This ends the proof of the Lemma.
\end{proof}

Now we  conclude the proof of Theorem \ref{soj}.

\begin{proof}

Since in the interior of $\xo,$ $\beta_0$ is a $C^\infty$ diffeomorphism between $C^\infty$ open manifolds,  $\beta_0^*\tilde\La$ is a $C^\infty$ Lagrangian manifold in the interior of $\mr_t\times \xo,$ and it is defined as
 \beq
\beta_0^*\tilde\La = \bigcup_{t_1\geq 0, t_2\geq 0} \exp t_2 \beta_0^*H_{Q_{ L}}\circ \exp t_1 \beta_0^*H_{Q_{R}} (\beta_0^*\Sigma),
 \eeq
where
\beq
\Sigma = \{(0, \tau, x,y, \xi,\eta, x, y,-\xi, -\eta) : x^2 \xi^2 + x^2 h(x, y, \eta) = \tau^2\}.
\eeq
  
In projective coordinates
 \begin{gather*}
 x', \;\ X= \frac{x}{x'}, \;\ Y=\frac{y-y'}{x'}, \;\ y',
 \end{gather*}
valid near $\ff$ and $L,$  $\beta_0^*\Sigma$
 can be written as
\beq
\beta_0^*\Sigma = \{(0,1,X, Y,\la,\mu, x', y' , \la', \mu'):  X = 1, Y = 0, \la' = \mu' = 0,  \la^2 + h(x', y', \mu) = \tau^2\},
\eeq
which is a $C^\infty$ submanifold of $T^*(\mr_t\times X\times_0 X)$ that extends smoothly up to the front face $\mr_t\times \ff=\{x'=0\}.$  Since $\beta_0^*\Sigma$ does not intersect either $\mr_s\times L$ or $\mr_s\times R,$ these properties do not change if we set $s=t+\log \rho_R+\log \rho_L,$  and hence
$\beta_1^*\Sigma$  is a $C^\infty$ submanifold of $T^*(\mbr_s \times \xo)$ that has a $C^\infty$ extension up to $\mr_s\times \ff.$   

 In the interior of $\mr_s\times \xo,$ $\beta_1^*Q_{L}$ and $\beta_1^*Q_R$ vanish on $\beta_1^*\tilde\La,$ and hence the integral curves of $H_{q_{L}}$ and 
 $H_{q_R}$ on $\beta_1^*\tilde \La$  coincide with the integral curves of $H_{\beta_1^*Q_{L}}$ and $H_{\beta_1^*Q_R}$ respectively. Therefore, in the interior of $\mr_s\times \xo$ and across to the front face,  $\beta_1^*\tilde\La$ is the union of integral curves of $H_{q_{L}}$ and $H_{q_R}$ emanating from $\beta_1^*\Sigma.$

  Since $Q_R$ and $Q_L$ do not depend on $t,$ it follows that $q_L$ and $q_R$ do not depend on $s,$ and hence $\sigma$ remains constant along the integral curves of $q_L$ and $q_R.$  Since $\sigma=\tau\not=0$ on $\beta_1^* \Sigma,$  it follows that $\sigma\not=0$ on  $\beta_1^*\tilde \La$ in the interior of $\mr_s \times \xo.$ However,  we have also shown that, up to the front face, in the region $\sigma=1,$ $H_{q_L}$ is transversal to $\mr_s\times L$ while $H_{q_R}$ is transversal  up to $\mr_s\times  R.$

  Recall from \eqref{liftedHA} that $H_{q_L}$ and $H_{q_R}$ are $C^\infty$ up to $\mr_s\times \ff$ and are tangent to $\mr_s\times \ff.$   So,  $\beta_1^* \tilde\La$ extends up to  $\mr_s \times \ff$ as the joint  flow-out of $\beta_1^*\Sigma$ by $H_{q_R}$ and $H_{q_L}.$

  So the integral curves of $H_{q_L}$ can be continued smoothly up to $\mr_s\times L$ and
the  integral curves of $H_{q_R}$ can be continued smoothly up to $\mr_s\times R.$   Therefore  $\beta_1^* \tilde{\Lambda}$ can be extended up to the face  $\{\rho_R=0\}$ because $H_{q_R}$ is tangent to $\beta_1^* \tilde{\Lambda}$ and transversal to $\{\rho_R=0\}.$  The same holds for the left face.   This shows that $\beta_1^* \tilde\La,$ which is in principle is defined in the interior of $\mr_s\times \xo,$ extends to a $C^\infty$ manifold up to $\p (\mr_s\times \xo) $ which intersects $\mr_s\times L$ and $\mr_s\times R$ transversally. See Fig. \ref{fig4}. 

We can make this more precise if we work suitable local symplectic coordinates valid near a point on the fiber over the corner $\ff\cap L \cap R.$   
We know that $R,$ $L$ and $\ff$ intersect transversally.  So one can choose local coordinates $x=(x_1,x_2,x_3, x')$ in $\mr^{2n+2}$ valid near  $\ff\cap L \cap R$ such that
\begin{gather*}
\ff=\{x_3=0\}, \;  R=\{x_1=0\}  \text{ and } L=\{x_2=0\}.
\end{gather*}
and that the symplectic form $\omega^0=d\sigma\wedge ds + d\xi \wedge d x.$  For example, this can be accomplished by using local coordinates defined in \eqref{eqc2} and setting $u=x_3,$ $w=x_2$ and $w'=x_1,$ $(y',Z)=x'.$ 

We know that $\beta_0^*\tilde\La$ is a Lagrangian submanifold of $T^*(\mr_s \times \mr^{2n+2})$ contained in $\{x_1> 0, \; x_2>0, \; x_3\geq 0\},$
 which intersects $\ff=\{x_3=0\}$ transversally.  There are commuting Hamilton vector fields $H_{q_R}$ and $H_{q_L}$ tangent to $\beta_0^*\tilde\La$ that are $C^\infty$ up to $\{x_1=0\}\cup\{x_2=0\}\cup \{x_3=0\},$  and as long as $\sigma\not=0,$ 
$H_{q_R}$ transversal to $R$ and tangent to $L$ and $\ff$ and $H_{q_L}$ is transversal to $L$ and tangent to $R$ and $\ff.$    Also, since $q_R$ and 
$q_L$ do not depend on $s,$ $\sigma$ remains constant along the integral curves of $H_{q_R}$ and $H_{q_L}.$

Let
\begin{gather*}
\mcf=T_{\{x_1=x_2=0\}}^*(\mr_s \times \{x: x_1>0, x_2>0, x_3\geq 0\}),
\end{gather*}
and let  $p=(s,\sigma, 0,\xi_1,0,\xi_2,x_3,\xi_3, x', \xi')),$ $\sigma\not=0,$  denote a point on  $\mcf,$ where

Since $q_R$ and $q_L$ do not depend on $s,$  $\sigma$ remains constant along the integral curves of  $H_{q_R}$ and $H_{q_L}.$  Moreover, in the region
$\sigma\not=0,$ the vector fields $H_{q_R}$ and $H_{q_L}$ are smooth, nondegenerate up to the boundaries.  $H_{q_R}$ is tangent to $\ff$ and $L,$
while $H_{q_L}$ is tangent to $\ff$ and $R.$ So, for $\eps$ small enough we  define
\begin{gather*}
\Psi_0: [0,\eps)\times [0,\eps) \times (\mcf \cap \{\sigma\not=0\} ) \longrightarrow U_0 \subset T^*(\mr_s\times \{ x_1\geq 0, x_2\geq 0, x_3\geq 0\})\\
\Psi_0(t_1,t_2,p) = \exp(-t_1 H_{q_R})\circ  \exp(-t_2 H_{q_L})(p),
\end{gather*}
and
\begin{gather*}
\Psi_1: [0,\eps)\times [0,\eps) \times (\mcf\cup \{\sigma\not=0\}) \longrightarrow U_1 \subset T^*(\mr_s\times \{ x_1\geq 0, x_2\geq 0, x_3\geq 0\})\\
\Psi_1(t_1,t_2,p) = \exp(-t_1 \p_{x_1})\circ  \exp(-t_2 \p_{x_2})(p),
\end{gather*}
Since the vector fields $H_{q_R},$ $H_{q_L}$ commute and  $\p_{x_1}$ and $\p_{x_2}$ commute, both maps are $C^\infty$ map and moreover, 
\begin{gather*}
\Psi_0^* H_{q_R}=-\p_{t_1}, \;\  \Psi_0^* H_{q_L}=-\p_{t_2}  \\
\Psi_1^* \p_{x_1}=-\p_{t_1}, \;\ \Psi_1^* \p_{x_2}=-\p_{t_2}  \\
\end{gather*}
Hence,
\begin{gather*}
\Psi=\Psi_0\circ \Psi_1^{-1}: U_1 \longrightarrow U_0, \\
\Psi^* H_{q_R}=-\p_{x_1},  \;\  \Psi^* H_{q_L}=-\p_{x_2}.
\end{gather*}
Moreover,  if $\omega^0$  is the symplectic form on $T^*(\mr \times \xo),$ in coordinates \eqref{eqc2} valid near $\mcf,$
\begin{gather*}
\Psi^*\omega^0= \omega^0.
\end{gather*}

Now $\Upsilon=\Psi^{-1}(\beta_0^*\tilde\La)$ is a $C^\infty$ Lagrangian in $\{ x_1> 0, \; x_2>0, x_3\geq 0\}$ which intersects $\{x_3=0\}$ transversally, and
both $\p_{x_1}$ and $\p_{x_2}$ are tangent to $\Upsilon.$  But this implies that  for any point $p \in \Upsilon,$ the integral curves of
 $\p_{x_j},$ $j=1,2$ starting at a point $p\in \Upsilon$ are contained in $\Upsilon.$  Therefore, for any $p=(x_1,\xi_1, x_2,\xi_2, x_3,\xi_3,x', \xi')\in \Upsilon,$ with 
 $x_1$ and $x_2$ small enough, the set
$\{x_1-t_1,\xi_1, x_2-t_2,\xi_2, x_3,\xi_3, x',\xi'\} \subset \Upsilon.$  By taking $t_1$ and $t_2$ large enough, this gives an extension $\overline{\Upsilon}$ of $\Upsilon$ to $\{x_1\leq 0\}\cup \{x_2\leq 0\}.$ Now $\Psi(\overline{\Upsilon})$ is the desired Lagrangian extension of $\beta_0^*\tilde\Lambda.$  Notice that in fact, it extends past the boundaries $\{x_1=0\}$ and $\{x_2=0\}.$  The construction in the other regions, away from the  co-dimension three corners follows by the same argument.

\epf

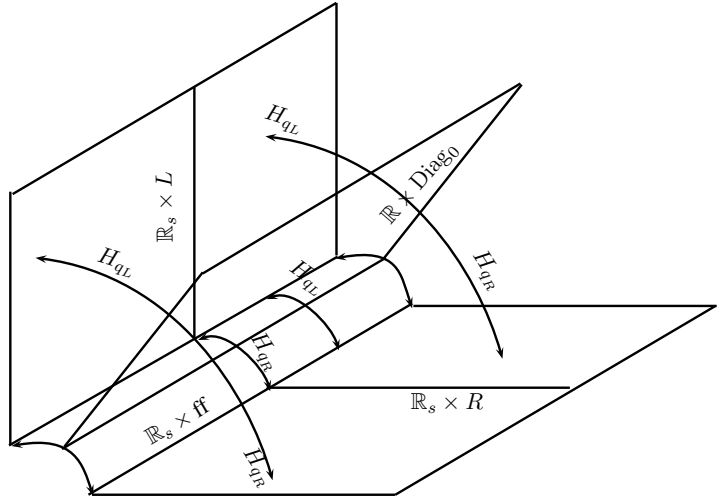
\begin{figure}
\scalebox{.8} 
{
\begin{pspicture}(0,-4.4649196)(11.78,4.4373207)
\psline[linewidth=0.04cm](3.08,3.0173206)(3.08,-1.1826793)
\psline[linewidth=0.04cm](4.28,-1.9826794)(9.32,-1.9826794)
\psline[linewidth=0.04cm](5.46,0.19732067)(0.0,-2.9426794)
\psline[linewidth=0.04cm](6.7,-0.58267933)(1.32,-3.7626793)
\psline[linewidth=0.04cm](1.38,-3.7626793)(6.42,-3.7626793)
\psline[linewidth=0.04cm](0.02170551,1.2773206)(0.02170551,-2.9226794)
\psline[linewidth=0.04cm](5.4417057,4.4173207)(5.4417057,0.21732067)
\psline[linewidth=0.04cm](6.72,-0.62267935)(11.76,-0.62267935)
\psline[linewidth=0.04cm](11.74,-0.62267935)(6.42,-3.7626793)
\psline[linewidth=0.04cm](0.06,1.2373207)(5.44,4.4173207)
\psline[linewidth=0.04cm](0.88,-3.0026793)(6.26,0.17732066)
\psline[linewidth=0.04cm](6.22,0.15732066)(8.52,3.0573206)
\psline[linewidth=0.04cm](0.92,-2.9626794)(3.22,-0.062679335)
\psline[linewidth=0.04cm](3.22,-0.08267934)(8.5,3.0773206)
\usefont{T1}{ptm}{m}{n}
\rput{49.219593}(3.4380233,-4.691262){\rput(6.84,1.4223206){$\mr\times \diag_0$}}
\usefont{T1}{ptm}{m}{n}
\rput(7.28,-2.2376792){$\mr_s\times R$}
\usefont{T1}{ptm}{m}{n}
\rput{91.17763}(3.7009497,-1.5710124){\rput(2.62,1.0423206){$\mr_s \times L$}}
\usefont{T1}{ptm}{m}{n}
\rput{30.563093}(-0.9177796,-1.7862711){\rput(2.81,-2.5576794){$\mr_s\times \ff$}}
\psbezier[linewidth=0.04,arrowsize=0.05291667cm 2.0,arrowlength=1.4,arrowinset=0.4]{<->}(4.26,2.1973207)(6.6,1.8973206)(7.8014956,-0.08267934)(8.2,-1.5026793)
\psbezier[linewidth=0.04,arrowsize=0.05291667cm 2.0,arrowlength=1.4,arrowinset=0.4]{<->}(0.44,0.17732066)(2.78,-0.12267934)(3.9814954,-2.1026793)(4.38,-3.5226793)
\psbezier[linewidth=0.04,arrowsize=0.05291667cm 2.0,arrowlength=1.4,arrowinset=0.4]{<->}(4.28,-0.50987935)(4.7561903,-0.32267934)(5.365714,-0.94667935)(5.48,-1.3626794)
\psbezier[linewidth=0.04,arrowsize=0.05291667cm 2.0,arrowlength=1.4,arrowinset=0.4]{<->}(3.14,-1.1898793)(3.6161904,-1.0026793)(4.225714,-1.6266793)(4.34,-2.0426793)
\psbezier[linewidth=0.04,arrowsize=0.05291667cm 2.0,arrowlength=1.4,arrowinset=0.4]{<->}(0.04,-2.9521713)(1.1,-2.7626793)(1.2316856,-3.2786305)(1.4,-3.7826793)
\psbezier[linewidth=0.04,arrowsize=0.05291667cm 2.0,arrowlength=1.4,arrowinset=0.4]{<->}(5.42,0.16782865)(6.48,0.35732067)(6.511686,-0.09863057)(6.7,-0.6626793)
\usefont{T1}{ptm}{m}{n}
\rput{-57.641018}(3.1470697,2.9342892){\rput(4.24,-1.3776793){$H_{q_{R}}$}}
\usefont{T1}{ptm}{m}{n}
\rput{-35.326515}(1.0112922,2.830469){\rput(4.95,-0.15767933){$H_{q_{L}}$}}
\usefont{T1}{ptm}{m}{n}
\rput{-17.133263}(-0.52405393,1.4557943){\rput(4.57,2.4823205){$H_{q_{L}}$}}
\usefont{T1}{ptm}{m}{n}
\rput{-57.641018}(3.71794,6.691819){\rput(7.94,-0.017679337){$H_{q_{R}}$}}
\usefont{T1}{ptm}{m}{n}
\rput{-17.133263}(0.05371211,0.53120047){\rput(1.79,0.10232066){$H_{q_{L}}$}}
\usefont{T1}{ptm}{m}{n}
\rput{-57.641018}(4.713142,1.9405507){\rput(4.12,-3.2976794){$H_{q_{R}}$}}
\end{pspicture} 
}
\caption{The integral curves of $H_{q_{R}}$ and $H_{q_{L}}$ in the region $\sigma=1.$}
\label{fig4}
\end{figure}
We remark that if $\Lambda_1$ is the extension of $\beta_1^* \tilde{\La},$ then, as mentioned in the introduction, the scattering relation is defined  as the intersection of $\Lambda_1$ to the corner $\mr_s\times \{\rho_R=\rho_L=0\}.$ Namely,
$SR=\Lambda_1\cap T_{(\mr_s \times\{\rho_R=\rho_L=0\})}^* (\mr_s\times \xo).$ Notice that the manifold $\{\rho_R=\rho_L=0\}=\p X \times_0 \p X$  is the blow-up of $\p X\times \p X$ induced by $\beta_0,$ and we let $\beta_{0,\p X}$ denote the blow-down map introduced in \cite{JS1}.

It is worth observing that Theorem \ref{soj} extends the result of \cite{SW} to the case of where  $(z,z') \rightarrow q \in \p X\times \p X.$    Let $(z_0,\zeta_0,z_0,-\zeta_0),$ $|\zeta_0|_{g^*(z_0)}=1,$  and let
\begin{gather*}
(z(t),\zeta(t), z'(t),\zeta'(t))= \exp( \frac{t}{2} H_{p_R}) \exp (\frac{t}{2} H_{p_L})( z_0,\zeta_0, z_0, -\zeta_0), \;\ t>0. 
\end{gather*}
 Let $\gamma$ be the curve defined by this equation joining $(z,\zeta)$ and $(z',\zeta').$  Suppose that  $\beta_0^{-1}(z(t), z'(t))  \rightarrow m \in \p X \times_0 \p X$  as $t \rightarrow \infty.$ The variable $s$ is a $C^\infty$ function on the Lagrangian $\beta_0^*( \tilde\Lambda)$ up to $\p(T^*(\mr_s \times \xo)).$  In particular  
$s(\beta_0^*(z,\zeta,z',\zeta'))$ is a $C^\infty$ function of the initial data $(z_0,\zeta_0),$ with $|\zeta_0|=1,$ up to the boundary.  In particular, its restriction to the boundary is a $C^\infty$ function, and we have
\begin{gather}
S(m,z_0,\zeta_0)=\lim_{t \rightarrow \infty} s\left(\beta_0^*((z(t),\zeta(t),z'(t), \zeta'(t))\right) \in C^\infty( \p X\times_0 \p X \times S^*\intx). \label{restric}
\end{gather}
In particular if there exist open subsets $U,V\subset \p X$ such that   $z(t)\rightarrow y\in U,$ $z'(t)\rightarrow  y'\in V,$ and $U\cap V=\emptyset,$  then $S(y,y',z_0,\zeta_0) \in C^\infty(U\times V \times S^*\intx).$

The definition of the sojourn time along a geodesic from \cite{SW} can be extended to the case where both points approach $\p X$ as
\begin{gather*}
S_{soj}=\beta_0^*(t + \log x +\log x')= s+2\log R,
\end{gather*}
and in particular $S_{soj}|_{\{\rho_R=\rho_L=0\}}= s|_{\{\rho_R=\rho_L=0\}} + 2\log R|_{\{\rho_R=\rho_L=0\}},$ and hence 
\begin{gather}
\left.\beta_0^* e^{-i\la S_{soj}}\right|_{\{\rho_R=\rho_L=0\}}= R^{-2i\la} H |_{\{\rho_R=\rho_L=0\}}, \;\ H\in C^\infty(\xo)\label{restric1}
\end{gather}
which according to \cite{JS1}  is a multiple the most singular term in the expansion of the Schwartz kernel  of  the scattering matrix.

\section{Asymptotics of the distance function}\label{distfunc}
 Let 
\begin{gather}
p(z,\zeta)=\ha(|\zeta|_{g^*(z)}^2-1), \text{ where } |\zeta|_{g^*(z)}^2=g(z)^*(\zeta,\zeta), \label{energy}
\end{gather}
and let 
\begin{gather}
S^*\diag=\{ (z,\zeta,z',\zeta') \in T^*(\intx \times \intx): \; z=z', \; \zeta=\-\zeta'; \;  p(z',\zeta')=0\}. \label{lag2N}
\end{gather}

We define the Lagrangian submanifold  of $T^*(\intx \times \intx),$
\begin{gather}
\Lambda= \bigcup_{t\geq 0} \exp(t H_{p}) S^*\diag=\{(z,\zeta,z',\zeta'): \exists \; t\geq 0, \; (z,\zeta)=\exp(t H_p)(z',\zeta'), \;\ p(z',\zeta')=0\} \label{Lag1}
\end{gather}

It is well-known that  integral curves of $H_p$ contained in $\{p=0\}$ project onto geodesics of the metric $g,$ see section 2.7 of \cite{AM}.  In other words, if $(z',\zeta')\in T^*\intx$ and $p(z',\zeta')=0,$  and if $\gamma(r)$ is a curve such that
\begin{gather*}
\frac{d}{dr} \gamma(r)=H_p(\gamma(r)), \\
\gamma(0)=(z',-\zeta').
\end{gather*}
Then,  with $\mcg$ given by \eqref{isomor},
\begin{gather}
\begin{gathered}
\gamma(r)= \left(\alpha(r), \mcg( \frac{d}{dr} \alpha(r)) \right), \text{ where } \\
\alpha(r)= \exp_{z'}(rv), \;\  v\in T_{z'} \intx, \;\ |v|_g^2=1, \;\ \mcg(v)= -\zeta',
\end{gathered}\label{geod-1}
\end{gather}
 $\exp_{z'}(r \bullet)$ denotes the exponential map on $T_{z'} \intx.$ 
  Now we assume that $(\intx,g)$ is geodesically convex.  In this case, the exponential map
$\exp_{z'}(r \bullet)$  is a global diffeomorphism for all $z'\in \intx$.   If 
\begin{gather*}
S_{z'}=\{v\in T_{z'} \intx:  |v|_g^2=1\},
\end{gather*}
the exponential map $\exp_{z'}$ gives geodesic normal coordinates about the point $z':$ 
\begin{gather*}
 S_{z'} \times [0,\infty) \longrightarrow \intx \\
(v,r) \mapsto  \exp_{z'} (r v).
\end{gather*}
The geodesic  is given by $\exp_{z'}(r v)=(z',v, r),$ and its tangent vector is $\p_r.$ Moreover, the metric $g$ takes the form
\begin{gather*}
g=dr^2+r^2 H(z,r, v, v),
\end{gather*}
and in this case $\mcg (\p_r)=dr,$ and according to \eqref{geod-1},  the integral curve of $H_p$ starting at $(z',-\zeta')$ is given by
\begin{gather}
\gamma(r)=(z',-\zeta', v, r, dr). \label{derr}
\end{gather}
So, the starting point of the curve $\gamma(r)$ is $(z',-\zeta')$ and its end point is $(z,\zeta)$ with $\zeta=dr.$  If one now reverses the role of $z$ and $z',$ and follows the geodesic in the opposite orientation, the starting point is $(z,-\zeta),$ and the end point is $(z',\zeta'),$ with $\zeta'=dr.$     
Then \eqref{Lag1} and \eqref{derr} give that

\begin{prop}\label{Lagrangian0} Let $(\intx,g)$ be a geodesically convex AHM. Then the manifold 
$\Lambda$ defined in \eqref{Lag1} is a $C^\infty$ embedded Lagrangian submanifold, and away from the diagonal $\Lambda$ is the graph of the differential of the distance function. In other words,
\begin{gather}\label{graph}
\Lambda \setminus S^*\diag=\{(z,\zeta,z',\zeta'):  \zeta'=-d_{z'}r(z,z'), \; \zeta= d_{z} r(z,z'), \text{ provided } (z,\zeta)\not= (z',-\zeta')\}.
\end{gather}
\end{prop}

\bpf[Proof of Theorem \ref{main}]  According to Proposition \ref{Lagrangian0}, away from the diagonal, the Lagrangian $\tilde \La$ defined in \eqref{eqlat}  satisfies
\beq
\tilde\La\cap \{\tau=1\}= \{(t, 1, z, \zeta, z', \zeta') \in T^*(\mbr_t\times\intx\times\intx):  (z,\zeta)= \exp(t H_p)(z',\zeta'), \;\ p(z',\zeta')=0\},
\eeq
On the other hand, in view of \eqref{graph},  away from the diagonal,
\begin{gather*}
\tilde\La\cap \{\tau=1\}= \{(t, \tau, z, \zeta, z', \zeta') \in T^*(\mbr_t\times\intx\times\intx):   \tau=1, \; t=r(z,z'),\; \zeta=d_z r(z,z'), \; \zeta'=d_{z'} r(z,z')\}.
\end{gather*}
Theorem \ref{soj},  guarantees that manifold $\tilde \La \cap \{\tau=1\}$ has a smooth extension to $T^*(\mbr_s\times\xo),$ where away from $\diag_0,$
\begin{equation}
s = t + \log \rho_L + \log \rho_R= r(z,z')+ \log \rho_R+\log \rho_L. \label{final}
\end{equation}
In particular, the right hand side of \eqref{final}  shows that  $s$ is in fact a function of the base variables only, and $s\in C^\infty(\xo\setminus \diag_0)$.  
\epf

%

\end{document}